\documentclass{amsart}
\usepackage{amssymb,amscd,verbatim, amsthm,amsmath,amsgen,xspace}

\usepackage[normalem]{ulem}

\usepackage{graphicx,color}

\newcommand\ep{{\epsilon}}

\newcommand\la{{\lambda}}

\def\im{\mathop{\hbox {Im}}\nolimits}

\def\ker{\mathop{\hbox{Ker}}\nolimits}

\def\mfrak{\mathfrak}

\newcommand{\pf}{\begin{proof}}
\newcommand{\epf}{\end{proof}}
\newcommand{\eq}{\begin{equation}}
\newcommand{\eeq}{\end{equation}}
\newcommand{\eqn}{\begin{equation*}}
\newcommand{\eeqn}{\end{equation*}}
\newcommand\bed{\begin{definition}}
\newcommand\ebed{\end{definition}}
\newcommand\bethm{\begin{theorem}}
\newcommand\ebethm{\end{theorem}}
\newcommand\bealigned{\begin{aligned}}
\newcommand\ebealigned{\end{aligned}}

\newcommand{\frg}{\mathfrak{g}}

\newcommand{\frk}{\mathfrak{k}}

\newcommand{\frp}{\mathfrak{p}}

\newcommand{\frz}{\mathfrak{z}}

\newcommand{\frsl}{\mathfrak{sl}}
\newcommand{\frsu}{\mathfrak{su}}

\newcommand{\bbC}{\mathbb{C}}

\newcommand{\bbN}{\mathbb{N}}

\newcommand{\bbZ}{\mathbb{Z}}

\newcommand{\Cas}{\operatorname{Cas}}
\newcommand{\half}{\frac{1}{2}}
\newcommand{\uqsl}{U_q(\frsl_2)}
\newcommand{\reduqsl}{U^{red}_q(\frsl_2)}
\newcommand{\uqslcp}{U_q(\frsl_2)\otimes C(\frp)}
\newcommand{\uqk}{U_q(\frk)}
\newcommand{\ki}{K^{-1}}
\newcommand{\Ki}{K^{-1}}
\newcommand{\qi}{q^{-1}}

\newtheorem{theorem}[equation]{Theorem}
\newtheorem{corollary}[equation]{Corollary}
\newtheorem{example}[equation]{Example}
\newtheorem{lemma}[equation]{Lemma}
\newtheorem{proposition}[equation]{Proposition}
\newtheorem{definition}[equation]{Definition}
\newtheorem{remark}[equation]{Remark}

\numberwithin{equation}{section}

\begin{document}
%\today

\bigskip
\title[Dirac operator for $\uqsl$]{Dirac operator and its cohomology for the quantum group $\uqsl$}

\author{Pavle Pand\v zi\'c}
	\address[Pand\v zi\'c]{Department of Mathematics, Faculty of Science, University of Zagreb, Bijeni\v cka 30, 10000 Zagreb,
          Croatia} 
	\email{pandzic@math.hr}
\thanks{P.~Pand\v zi\'c was supported by grant no. 4176 of the
  Croatian Science Foundation and by the QuantiXLie Center of Excellence.} 
\author{Petr Somberg}
\address[Somberg]{Mathematical Institute MFF UK\\
Sokolovsk\'a 83, 18000 Praha 8 - Karl\'{\i}n, Czech Republic}
\email{somberg@karlin.mff.cuni.cz}
\thanks{P.~Somberg was supported by grant GA \v CR P201/12/G028.}
\begin{abstract} We introduce a Dirac operator $D$ for the quantum group $\uqsl$, as an element of the tensor product of $\uqsl$ with the Clifford algebra on two generators. We study the properties of $D$, including an analogue of Vogan's conjecture. We compute the cohomology of $D$ acting on various $\uqsl$-modules.
\end{abstract}


\maketitle

%%%%%%%%%%%%%%%%%%%%%%%%%%%%%%%%%%%%%%%%%%%%%%%%%%%%%%%%%%%%%%%%%%%%%%%%%%%%%%%%%%%%%%%%%%%
\section{Introduction}\label{section1}

 Operators of Dirac type were introduced into the representation theory of 
classical reductive Lie groups by Partasarathy \cite{Par} in order to construct the discrete 
series representations as well as to study the unitary $({\mathfrak g},K)$-modules. 
In 1997 D. Vogan \cite{V}
introduced the notion of Dirac cohomology for a wide class of modules with an action of a Dirac operator. He conjectured a relationship between 
$\tilde{K}$-types constituting the Dirac 
cohomology of a $({\mathfrak g},K)$-module $M$ and the infinitesimal character of $M$. This conjecture was proved by Huang and Pand\v zi\'c, \cite{HP}, \cite{HP2}. 
The framework was consequently extended and analogous results have been proved in several other settings:
\begin{itemize} 
\item For a reductive Lie algebra $\frg$ and an arbitrary quadratic subalgebra \cite{Ko};
\item For a basic classical Lie superalgebra \cite{HP3};
\item For an affine Lie algebra \cite{KMP};
\item For noncommutative equivariant cohomology \cite{AlMe}, \cite{Ku};
\item For graded affine Hecke algebras, with applications to p-adic groups \cite{BCT};
\item For symplectic reflection algebras \cite{C}.
\end{itemize}

Quantum groups as a mathematical structure have their origin in many 
problems studied in theoretical physics, e.g. the solution of Yang-Baxter 
equation, the description of monodromy of the vertex operators in 
conformal field theory, the integrable systems, etc.
They naturally arise as Hopf algebras depending on an auxiliary 
parameter $q$ (or $h$), which specialize to the universal enveloping algebras 
of certain Lie algebras (quite often semisimple) for $q=1$ 
(or $h=0$.) 
As in the case of classical Lie groups, appearance of quantum groups as the source of 
symmetries comes through the notion of their representations. The
central problem is then the description and characterization of representations (modules) of a quantum group. 
 
In the present article, we initiate and study Dirac operators and cohomology for quantum groups. Namely, we introduce a Dirac operator and its cohomology theory
in the case of the quantized universal enveloping 
algebra $\uqsl$ associated to the simple Lie algebra ${\mathfrak sl}(2,{\mathbb C})$, 
and prove various
structural results, both when $q$ is not a root of unity, and when $q$ is a root of unity. We note that Dirac operators in a different but related setting of quantum $SU(2)$ group and quantum sphere were studied in \cite{BK}.

It is more difficult to obtain this theory for general quantized enveloping algebras $U_q({\mathfrak g})$.
We plan to do it in near future, using the concept of braided Lie algebras and braided Killing forms \cite{M}, as well as a version of the analogues of Harish-Chandra modules studied in \cite{Le}. 
One of the problems to overcome is to find a suitable analogue for the classical center of the enveloping algebra.

We now briefly review the content of our article. In Section \ref{section2} we
first recall the definition and basic structural properties of the Hopf algebra $\uqsl$. As an algebra, it is generated by $K,\ki,E,F$ satisfying appropriate commutation relations.

On the other hand, we also consider the classical Lie algebra $\frg=\frsl(2,\bbC)$ with the usual basis $H,E_0,F_0$. 
On the subspace $\frp$ of $\frg$ spanned by $E_0$ and $F_0$, we consider the trace form $B$, so that $E_0$ and $F_0$ are isotropic and dual to each other. To this quadratic space we associate the Clifford algebra $C(\frp)$. To avoid notational confusion, we will denote by $c$ the embedding of $\frp$ (which is a subspace of $\frg$ into $C(\frp)$. Then the standard generators of $C(\frp)$ will be denoted by $c(E_0)$ and $c(F_0)$. 

The Dirac operator is defined as
\[
D=E\otimes c(F_0)+F\otimes c(E_0) \qquad\quad \in \uqsl\otimes C(B).
\]
We compute a Parthasarathy
type formula for $D^2$ and discuss the analogue of Vogan's conjecture. In Section \ref{section3} we highlight the
relationship of our Dirac operator to the representation theory of $\uqsl$,
and introduce the notion of Dirac cohomology for any $\uqsl$-module. 
In Section \ref{section4} we review the classification of finite-dimensional
$\uqsl$-modules in the case when $q$ is not a root of unity, and 
determine the Dirac cohomology for all these representations. 
In Section \ref{section5} we recall the classification of (necessarily 
finite-dimensional) irreducible $\uqsl$-modules in the case when $q$ 
is a root of unity, and as in the previous section 
compute their Dirac cohomology. In Section \ref{section6} 
we consider some further examples, including analogues of Verma modules and some finite-dimensional indecomposable reducible modules (for $q$ a root of unity). 

Throughout the article we denote by ${\mathbb N}$ the set of 
integers $\{1,2,\dots\}$, and by ${\mathbb N}_0$ the integers together 
with zero. The symbol $[n]=[n]_q$ is the $q$-integer, defined by 
\[
[n]_q=\frac{q^n-q^{-n}}{q-q^{-1}},\qquad n\in{\mathbb N}_0.
\] 
A complex vector space spanned by $v_1,\ldots ,v_k$ is denoted by $\langle v_1,\ldots ,v_k\rangle$, except if $k=1$, when we denote the space by $\bbC v_1$.

%%%%%%%%%%%%%%%%%%%%%%%%%%%%%%%%%%%%%%%%%%%%%%%%%%%%%%%%%%%%%%%%%%%%%%%%%%%%%%%%%%%%%%%%%%%%

\section{Dirac operator for $\uqsl$}\label{section2}

We follow the conventions and notation from \cite{KS}; another good introduction to quantum groups is \cite{L2}. Let $q$ be a fixed complex number not equal to $0$ or $\pm 1$. Let $\uqsl=U_q(\mfrak{sl}(2,{\mathbb C}))$ 
be the associative unital algebra over $\bbC$ generated by 
\[
K, K^{-1}, E, F,
\] 
with relations 
\begin{eqnarray}
\label{defuq}
& KK^{-1}=1=K^{-1}K; \nonumber \\
& KE=q^2EK,\qquad KF=q^{-2}FK; \\ 
& EF-FE=\frac{K-K^{-1}}{q-q^{-1}}\nonumber.
\end{eqnarray}
We also consider the classical Lie algebra $\frg=\frsl(2,\bbC)$, with basis $H,E_0,F_0$ satisfying the usual commutation relations
\[
[H,E_0]=2E_0,\quad [H,F_0]=-2F_0,\quad [E_0,F_0]=H.
\]
As $q$ varies, the family of algebras $\uqsl$ can be thought of as a deformation of the universal enveloping algebra $U(\frg)$ of $\frg$, and $U(\frg)$ is the ``classical limit" of $\uqsl$ for $q\to 1$. Loosely speaking, one can connect $K$ and $H$ by 
\[
K=q^H=e^{hH},
\]
where $h$ is related to $q$ by $q=e^h$ ($q\to 1$ is equivalent to $h\to 0$). On the other hand, passing to the classical limit takes $E$ to $E_0$ and $F$ to $F_0$.  

With the following definitions of coproduct, counit and antipode, $\uqsl$ becomes a Hopf algebra. The coproduct is the algebra homomorphism
\begin{align*}
\triangle:\, \uqsl\longrightarrow 
\uqsl\otimes \uqsl
\end{align*}
given on generators by
\begin{eqnarray} 
&\triangle(K^{\pm 1})=K^{\pm 1}\otimes K^{\pm 1},\nonumber\\ 
&\triangle(E)=E\otimes K+1\otimes E,\\ 
&\triangle(F)=F\otimes 1+K^{-1}\otimes\nonumber F.
\end{eqnarray}
The counit is the algebra homomorphism $\ep:\uqsl\to\bbC$ given on generators by
\begin{eqnarray}
\ep(K)=\ep(K^{-1})=1,\qquad \ep(E)=\ep(F)=0.
\end{eqnarray}
The antipode $S$ is the antiautomorphism of the algebra $\uqsl$ given on generators by
\begin{eqnarray}
&S(K)=K^{-1},\quad S(K^{-1})=K,\\
&S(E)=-EK^{-1},\quad S(F)=-KF.\nonumber
\end{eqnarray}

If $q$ is not a root of unity, the center $\frz$ of $\uqsl$ is generated by the Casimir element 
\begin{align}\label{Casimir'}
\mathrm{Cas}'_q=EF+\frac{qK^{-1}+q^{-1}K}{(q-q^{-1})^2}.
\end{align}
However, we prefer to use the following normalized form of $\Cas_q'$: 
\begin{multline}\label{cassl2q}
\qquad\qquad\mathrm{Cas}_q=2EF+\frac{2qK^{-1}+2q^{-1}K-2(q+\qi)}{(q-q^{-1})^2}\\ =EF+FE+\frac{(q+q^{-1})(K+K^{-1}-2)}{(q-q^{-1})^2}\qquad\qquad.
\end{multline}
It is easy to check, by expanding $K=e^{hH}$, $K^{-1}=e^{-hH}$, $q=e^h$ into Taylor series, that  
\[
\lim\limits_{q\rightarrow 1}\mathrm{Cas}_q=\mathrm{Cas},
\]
where $\mathrm{Cas}=E_0F_0+F_0E_0+\half H^2$ is the Casimir element for $\mfrak{sl}(2,{\mathbb C})$ with respect to the trace form. On the other hand, $\Cas'_q$ has no limit for $q\to 1$.

If $q$ is a root of unity, the center of $\uqsl$ is much bigger, as described in the following proposition \cite{DKP}, \cite{KS}, 3.3.1., Proposition 15.

\begin{proposition}
\label{center}
(i) If $q$ is not a root of unity, then the center of $\uqsl$ consists of polynomials in $\Cas_q$.

(ii) Suppose $q$ is a primitive $p'$-th root of unity, with $p'\geq 3$. Set
\[
p=p',\quad p'\text{ odd};\qquad p=p'/2,\quad p'\text{ even}.
\]
Then the center of $\uqsl$ is generated by $E^p, F^p, K^p, K^{-p}$ and $\Cas_q$. \qed
\end{proposition}

Let $\frg=\frk\oplus\frp$ be the Cartan decomposition of $\frg=\frsl(2,\bbC)$ corresponding to the real form $\frsu(1,1)$; so $\frk=\bbC H$, and $\frp$ is spanned by $E_0$ and $F_0$. Let $C(\frp)$ be the Clifford algebra of $\frp$ with respect to the trace form $B$. It is the associative unital algebra with generators $c(E_0)$ and $c(F_0)$ corresponding to $E_0$ and $F_0$ under the embedding $c:\frp\to C(\frp)$, with relations
\begin{align} 
\label{cliffalg}
 c(E_0)^2=0,\quad c(F_0)^2=0,\quad c(E_0)c(F_0)+c(F_0)c(E_0)=1.
\end{align}
The following lemma is easy to check. It is also 
well known in much greater generality.

\begin{lemma}
\label{emb cliff}
The map $\alpha:\frk\to C(\frp)$ defined by $\alpha(H)=2c(E_0)c(F_0)-1$ is a homomorphism of Lie algebras. Furthermore, for any $X\in\frk$ and $Y\in\frp$ we have
\[
c[X,Y]=[\alpha(X),c(Y)],
\]
with the bracket on the left hand side computed in $\frg$, and the bracket on the right hand side computed in $C(\frp)$.
\qed
\end{lemma} 

The reason why we are using the ordinary and not ``quantized" Clifford algebra is the fact that the Clifford algebra is not deformable in the category of associative unital algebras; see \cite{MPU}. Indeed, the following computation shows that
we can define $\alpha(K)=q^{\alpha(H)}=e^{h\alpha(H)}$ as an element of $C(\frp)$. In contrast, $K$ can not be found in the classical enveloping algebra $U(\frg)$. 

To compute $e^{h\alpha(H)}$ in $C(\frp)$, we first note that 
\[
(c(E_0)c(F_0))^2=c(E_0)c(F_0)c(E_0)c(F_0)=c(E_0)(-c(E_0)c(F_0)+1)c(F_0)=c(E_0)c(F_0),
\]
and therefore
\[
\alpha(H)^2=(2c(E_0)c(F_0)-1)^2=4(c(E_0)c(F_0))^2-4c(E_0)c(F_0)+1=1.
\]
It follows that for any integer $n\geq 0$, $\alpha(H)^{2n}=1$ and $\alpha(H)^{2n+1}=\alpha(H)$, and hence
\begin{multline*}
\alpha(K)=e^{h\alpha(H)}=\sum_{n=0}^\infty \frac{1}{n!}h^n\alpha(H)^n = \cosh h + \sinh h\,\, \alpha(H)=\\ (\cosh h-\sinh h) + 2\cosh h\,\, c(E_0)c(F_0) = \qi+(q-\qi)c(E_0)c(F_0). 
\end{multline*}
Similarly,
\[
\alpha(K^{-1})=e^{-h\alpha(H)}=q-(q-q^{-1})c(E_0)c(F_0).
\]
This and a little more computation immediately implies the following proposition.

\begin{proposition}
\label{qcliff}
The elements 
\[
\alpha(K)=\qi+(q-\qi)c(E_0)c(F_0),\qquad \alpha(\ki)=q-(q-\qi)c(E_0)c(F_0)
\]
of $C(\frp)$ define an algebra homomorphism $\alpha$ from the subalgebra $U_q(\frk)=\bbC[K,\ki]$ of $\uqsl$ into $C(\frp)$. Furthermore, the elements $\alpha(K)$ and $\alpha(\ki)=\alpha(K)^{-1}$ of $C(\frp)$ satisfy the relations
\begin{eqnarray}
\label{Kconjcliff}
& \alpha(K)c(E_0)\alpha(K)^{-1} =q^2c(E_0);\quad \alpha(K)c(F_0)\alpha(K)^{-1} =q^{-2}c(F_0);\\
&\frac{\alpha(K)-\alpha(K)^{-1}}{q-q^{-1}}=2c(E_0)c(F_0)-1=\alpha(H). \nonumber
\end{eqnarray}
\qed
\end{proposition}

We note that there is a family of graded coproducts on $C(\frp)$ introduced in \cite{Pan}, but the homomorphism $\alpha$ does not have any good property with respect to coproducts.

\begin{corollary}
\label{diag emb}
Let $\alpha:U_q(\frk)\to C(\frp)$ be as in Proposition \ref{qcliff}, and let $i:U_q(\frk)\to\uqsl$ be the inclusion map. Then 
\[
\delta=(i\otimes \alpha)\circ \triangle: U_q(\frk)\to \uqsl\otimes C(\frp)
\]
is an injective homomorphism of algebras, which we call the diagonal embedding. Explicitly, $\delta$ is given by
\begin{eqnarray}
\label{deltaK}
& \delta(K)=q^{-1}K\otimes 1+(q-q^{-1})K\otimes c(E_0)c(F_0), \\
& \delta(K^{-1})=qK^{-1}\otimes 1-(q-q^{-1})K^{-1}\otimes c(E_0)c(F_0). \nonumber
\end{eqnarray}
\qed
\end{corollary}

We introduce the following Casimir element for $U_q(\mathfrak{k})$:
\begin{align} 
\mathrm{Cas}_q(\mathfrak{k})= 
\frac{(q+q^{-1})(K+K^{-1}-2)}{(q-q^{-1})^2}.
\end{align}
Then (\ref{cassl2q}) implies
\[
\Cas_q=\Cas_q(\frk)+EF+FE.
\]
A short computation using (\ref{deltaK}) gives
\begin{equation}
\label{cas delta k}
\delta(\Cas_q(\frk))=\frac{(q+q^{-1})(q^{-1}K+qK^{-1}-2)}{(q-q^{-1})^2}\otimes 1
 +(q+q^{-1})\frac{K-K^{-1}}{q-q^{-1}}\otimes c(E_0)c(F_0).
\end{equation}

\begin{definition}
\label{def dirac}
The Dirac operator $D$ for $\uqsl$ is
\eqn 
D=E\otimes c(F_0)+F\otimes c(E_0) \qquad
\in \uqsl\otimes C(\mathfrak{p}).
\eeqn
\end{definition}

The following lemma is proved by an easy computation.

\begin{lemma}
\label{D k inv}
The Dirac operator $D$ is invariant under the action of 
$\uqk=\bbC[K,\ki]$ on $\uqslcp$, defined by letting $K$ and $\ki$ act by conjugation in the first factor and conjugation by $\alpha(K)$, $\alpha(\ki)$ in the second factor, as in (\ref{defuq}) and (\ref{Kconjcliff}).
\end{lemma}

The following formula for $D^2$ is an analogue of the well known formula of Parthasarathy \cite{Par}. 

\begin{proposition}
\label{D squared}
The square of the Dirac operator $D$ is
\[
D^2=\half\Cas_q\otimes 1 -\frac{1}{q+\qi}\delta(\Cas_q(\frk))+\frac{q+\qi-2}{(q-\qi)^2}\otimes 1.
\]
\end{proposition}
\pf
We can write $D^2$ as
\[
D^2 =(E\otimes c(F_0)+F\otimes c(E_0))^2=-(EF-FE)\otimes c(E_0)c(F_0)+EF\otimes 1.
\]
By the defining relations (\ref{defuq}) for $\uqsl$, and by (\ref{cas delta k}), we can write
\begin{multline}
-(EF-FE)\otimes c(E_0)c(F_0)=-\frac{(K-K^{-1})}{(q-q^{-1})}\otimes c(E_0)c(F_0) 
\\=
-\frac{1}{q+\qi}\delta(\Cas_q(\frk))+\frac{\qi K+q\Ki-2}{(q-\qi)^2}\otimes 1.
\end{multline}
On the other hand, using (\ref{cassl2q}) we can write
\[
EF\otimes 1=\left(\half\Cas_q-\frac{q\Ki+\qi K-(q+\qi)}{(q-\qi)^2}\right)\otimes 1.
\]
Adding up these two equalities proves the proposition.
\epf

\begin{remark}
{\rm
It is easy to check that the classical limit of the above formula for $D^2$ gives
\[
\lim_{q\to 1} D^2=\half\Cas-\half\delta(\Cas(\frk))+\frac{1}{4},
\]
where $\delta:U(\frk)\to U(\frg)\otimes C(\frp)$ is the diagonal embedding in the classical setting.
This is compatible with the classical formula of Partasarathy \cite{Par}. See \cite{HP2} for more details about the classical setting; note that the formulas are not completely equal because here we are using different conventions for $C(\frp)$ and its action on the spin module $S$.
}
\end{remark}

Recall that the center $\frz$ of $\uqsl$ was described in Proposition \eqref{center}; note that $\uqk$ is abelian and hence equal to its center.
We will use a $q$-analogue of the classical Harish-Chandra homomorphism introduced in \cite{DKP}, \S 2. By the Poincar\'e-Birkhoff-Witt Theorem (\cite{KS}, 3.1.1., Proposition 1, and 6.1.5, Theorem 14), the monomials
\eq
\label{pbw}
E^iK^jF^k,\qquad i,k\in\bbN_0,\, j\in\bbZ
\eeq
form a vector space basis for $\uqsl$. It follows that
\[
\uqsl=\uqk\oplus(E\uqsl + \uqsl F),
\]
and we define $\mu:\uqsl\to\uqk$ to be the projection along $E\uqsl + \uqsl F$. Since any element of $\frz$ commutes with $K$, its expression in the basis \eqref{pbw} only contains monomials with $i=k$. In particular, 
\[
\frz\cap (E\uqsl + \uqsl F) =\frz\cap E\uqsl =\frz\cap \uqsl F=\frz\cap E\uqsl F.
\]
This implies that $\frz\cap (E\uqsl + \uqsl F)$ is an ideal in $\frz$, and hence $\mu:\frz\to\uqk$ is an algebra homomorphism. Now we compose $\mu$ with the algebra automorphism 
\[
\sigma:\uqk\to\uqk,\qquad \sigma(K^j)=q^jK^j,
\]
which is a multiplicative version of the classical $\rho$-shift. We get 
\[
\gamma=\sigma\circ\mu:\frz\to\uqk
\]
analogous to the Harish-Chandra homomorphism.

Composing with the diagonal embedding $\delta$, we get an algebra homomorphism
\eq
\label{zeta}
\zeta=\delta\circ\gamma:\frz\to\delta(\uqk).
\eeq
The following theorem is an analogue of Vogan's conjecture \cite{V}, \cite{HP}, for the quantum group $\uqsl$. 

\begin{theorem}
\label{vogconj}
Let $\zeta:\frz\to \delta(\uqk)$ be the algebra homomorphism given by \eqref{zeta}. Then for any $z\in \frz$, there is a $\uqk$-invariant element $a\in \uqsl \otimes C(\frp)$ such that
\begin{align}\label{homotcond}
z \otimes 1 = \zeta(z) + Da + aD .
\end{align}
\end{theorem}
\pf 
If we have 
\begin{align} 
z_1 \otimes 1 = \zeta(z_1) + Da_1 + a_1D \quad\text{and}\quad
z_2 \otimes 1 = \zeta(z_2) + Da_2 + a_2D, 
\end{align}
then we also have
\begin{align}\label{multiplicstructure}
z_1z_2 \otimes 1 = \zeta(z_1)\zeta(z_2) 
+ D(\zeta(z_1)a_2+\zeta(z_2)a_1+a_1a_2)
\nonumber \\
+(\zeta(z_1)a_2+\zeta(z_2)a_1+a_1a_2)D. 
\end{align}
Since $\zeta$ is a homomorphism, we see that if the claim of the theorem is true for $z_1$ and $z_2$, it is also true for $z_1z_2$. So it is sufficient to prove the claim for generators of $\frz$, which are given in Proposition \ref{center}.

If $z=\Cas_q$ is the Casimir element, then Proposition \ref{D squared} implies
\[
\Cas_q\otimes 1 =\frac{2}{q+\qi}\delta(\Cas_q(\frk))-2\frac{q+\qi-2}{(q-\qi)^2}\otimes 1+2D^2.
\]
A short computation shows that 
\[
\zeta(\Cas_q)=\frac{2}{q+\qi}\delta(\Cas_q(\frk))-2\frac{q+\qi-2}{(q-\qi)^2}\otimes 1,
\]
so \eqref{homotcond} is true with $a=D$ (which is $\uqk$-invariant by Lemma \ref{D k inv}). This proves the theorem in case $q$ is not a root of unity.

Assume now that $q$ is a primitive $p'$-th root of unity, and let as before $p=p'$ if $p'$ is odd, and $p=p'/2$ if $p'$ is even. 
By Proposition \ref{center}, $\frz$ is generated by $\Cas_q$,  $E^p, F^p, K^p,$ and  $K^{-p}$. We already proved \eqref{homotcond} for $\Cas_q$, and now we have to prove it for $z$ equal to $E^p, F^p, K^p,$ or  $K^{-p}$. For $z=E^p$, we have 
\begin{align} 
E^p\otimes 1 = D(E^{p-1}\otimes c(E_0)) + (E^{p-1}\otimes c(E_0))D,
\end{align}
and since $\zeta(E^p)=0$, we see that \eqref{homotcond} holds with $a=E^{p-1}\otimes c(E_0)$ (which is $\uqk$-invariant sice $q^{2p}=1$). Similarly, \eqref{homotcond} holds
 for $z=F^p$ with $a=F^{p-1}\otimes c(F_0)$.

For $z=K^p$, we use \eqref{cliffalg} to compute 
\begin{align}
\zeta(K^p)=\delta(q^pK^p)=q^p\delta(K)^p=q^p\sum\limits_{i=0}^p\binom{p}{i}q^{-i}K^i(q-q^{-1})^{p-i}
K^{p-i}\otimes (c(E_0)c(F_0))^{p-i}
\nonumber \\
=K^p\otimes 1 + q^p\Big(\sum\limits_{i=0}^{p-1}\binom{p}{i}q^{-i}(q-q^{-1})^{p-i}\Big)
(K^{p}\otimes c(E_0)c(F_0))
\nonumber \\ 
= K^p\otimes 1 + 
q^p\Big(\sum\limits_{i=0}^{p}\binom{p}{i}q^{-i}(q-q^{-1})^{p-i}\Big)
(K^{p}\otimes c(E_0)c(F_0)) 
- K^p\otimes c(E_0)c(F_0) 
\nonumber \\
=K^p\otimes 1 + q^p(\qi+(q-\qi))^p (K^{p}\otimes c(E_0)c(F_0)) 
- K^p\otimes c(E_0)c(F_0) 
= K^p\otimes 1
\end{align}
due to $q^{2p}=q^{p'}=1$.
Hence \eqref{homotcond} holds with $a=0$. The case of $K^{-p}$ is similar.
\epf

%%%%%%%%%%%%%%%%%%%%%%%%%%%%%%%%%%%%%%%%%%%%%%%%%%%%%%%%%%%%%%%%%%%%%%%
\section{Dirac cohomology for $\uqsl$-modules}\label{section3}

Let $M$ be a module for $\uqsl$. In order to make the Dirac operator $D$ act, we should tensor $M$ with a module for $C(\frp)$; it is however well known (see for example \cite{HP2}, Chapter 2), that the only simple $C(\frp)$-module is the spin module $S$, and that all other modules are direct sums of copies of $S$. So the natural $C(\frp)$-module to use is $S$. Recall that 
\[
S=\bbC s_{-1}\oplus \bbC s_1,
\]
with $c(E_0)$ and $c(F_0)$ acting by
\begin{align*}
 & c(E_0)s_{-1} =s_{1},\, c(E_0)s_{1} =0,\\
 & c(F_0)s_{-1} =0,\, c(F_0)s_{1} =s_{-1}.
\end{align*}
It follows from Proposition \ref{qcliff} that $S$ is also a module over $U_q(\frk)=\bbC[K,\ki]$, through the homomorphism $\alpha:U_q(\frk)\to C(\frp)$. In other words, $K$ and $K^{-1}$ act by 
\begin{align}
Ks_{-1} =q^{-1}s_{-1},\quad Ks_{1} =qs_1,\quad K^{-1}s_{-1} =qs_{-1},\quad K^{-1}s_{1} =q^{-1}s_{1}.
\end{align}

The algebra $\uqsl\otimes C(\frp)$ acts on $M\otimes S$ by
\[
(u\otimes a)(m\otimes s)=um\otimes as,\qquad u\in\uqsl,\ a\in C(\frp),\ m\in M,\ s\in S.
\]
It follows that $D\in\uqsl\otimes C(\frp)$ acts on $M\otimes S$, and we define the Dirac cohomology of $M$ to be the vector space
\eq
\label{dircoh}
H_D(M)=\ker (D) \big/ (\im (D)\cap\ker (D)).
\eeq
Since $D$ is invariant under $\uqk$ by Lemma \ref{D k inv}, 
$H_D(M)$ is a $\uqk$-module.

We have an $\uqk$-module isomorphism 
\begin{align}
M\otimes S\cong M\otimes s_{-1}\oplus M\otimes s_1,
\end{align}
and the Dirac operator $D$ acts on $M\otimes S$ by
\begin{align}
D(m\otimes s_{-1})=Fm\otimes s_1,\quad D(m\otimes s_{1})=Em\otimes s_{-1}, \quad \quad m\in M .
\end{align}
It follows that 
\begin{align}
& \ker(D)=\ker(F)\otimes s_{-1} \oplus \ker(E)\otimes s_{1},
\nonumber \\
& \im(D)=\im(F)\otimes s_1\oplus \im(E)\otimes s_{-1},
\end{align}
and hence
\begin{align}\label{Diraccoh}
H_D(M) & = \ker(F)/(\ker(F)\cap\im(E))\otimes s_{-1}
\oplus \ker(E)/(\ker(E)\cap\im(F))\otimes s_{1}.
\end{align}

As before, let $\frz$ be the center of $\uqsl$ and let $\chi:\frz\to\bbC$ be a character (i.e., $\chi$ is an algebra homomorphism). We say that a $\uqsl$-module $M$ has infinitesimal character $\chi$, if every $z\in\frz$ acts on $M$ by the scalar $\chi(z)$. 

On the other hand, on all modules we will consider in this paper, the commutative algebra $\uqk=\bbC[K,\ki]$ acts semisimply, with eigenvalues given by integer powers of $q$. In other words, we consider $\uqsl$-modules 
\eq
\label{Kss}
M=\bigoplus_{i\in\bbZ} M(i),
\eeq
with $\uqk$ acting on $M(i)$ by the character $\chi^\frk_i$, $i\in\bbZ$, where
\[
\chi^\frk_i(K)=q^i.
\]
Then we have the following standard consequence of the analogue of Vogan's conjecture given by Theorem \ref{vogconj}.

\begin{theorem}
\label{vogconjmod}
Let $M$ be a $\uqsl$ module as in \eqref{Kss}, and assume $M$ has infinitesimal character $\chi$. Assume that the Dirac cohomology $H_D(M)$ has a $\uqk$-submodule $E_i$ on which $\uqk$ acts by the character $\chi^\frk_i$. Then $\chi=\chi^\frk_i\circ\zeta$, for $\zeta:\frz\to\uqk$ defined above Theorem \ref{vogconj}.
\end{theorem}
\pf By Theorem \eqref{vogconj}, for any $z\in\frz$, we have
\eq
\label{dircohpf}
z\otimes 1=\zeta(z)+Da+aD.
\eeq
Let $v\in\ker(D)\subset M\otimes S$ be a representative of a
nonzero element of $E_i$. Then $Dv=0$ and so \eqref{dircohpf} implies
\[
\chi(z)v=\chi^\frk_i(\zeta(z))v+Dav
\]
Since the class of $v$ in $H_D(M)$ is nonzero, $v\notin\im(D)$. It follows that $\chi(z)=\chi^\frk_i(\zeta(z))$. This implies the claim.
\epf 
%%%%%%%%%%%%%%%%%%%%%%%%%%%%%%%%%%%%%%%%%%%%%%%%%%%%%%%%%%%%%%%%%%%%%%%%%%%
\section{Irreducible finite-dimensional $\uqsl$-modules for $q$ not a root of unity}\label{modnrootof1}
\label{section4}
Throughout this section we assume $q$ is not a root of unity.
We recall the classification of irreducible finite-dimensional $U_q(\mathfrak{sl}_2)$-modules as described in \cite{KS}, Section 3.2. We use the notation $(T,V)$ for the representation $T$ on the complex vector space $V$. 

Let $k\in\half\bbN_0$ and let $\omega\in\{1,-1\}$. Let $V_{k}$ be a $(2k+1)$-dimensional 
complex vector space,
\[
V_k=\langle v_{-k},v_{-k+1},\dots,v_{k-1},v_k \rangle;
\]
set also $v_{-k-1}=v_{k+1}=0$. 
The representation $(T_{\omega,k},V_k)$ of $U_q(\mathfrak{sl}_2)$ is defined by setting
\begin{align}
\label{def vlo}
& T_{\omega,k}(K)v_m=\omega q^{2m}v_m,\nonumber \\
& T_{\omega,k}(E)v_m=\sqrt{[k-m][k+m+1]}v_{m+1}, \nonumber \\
& T_{\omega,k}(F)v_m=\omega\sqrt{[k+m][k-m+1]}v_{m-1}
\end{align}
for all $m\in\{-k,\ldots,k\}$. The square root is chosen to have argument in $[0,\pi\rangle$. Recall that 
\[
[n]=[n]_q=\frac{q^n-q^{-n}}{q-\qi}.
\]
\begin{theorem}
\label{cl fd nonroot} With notation as above,

(i) For all $(\omega,k)\in\{1,-1\}\times \half\bbN_0$, $T_{\omega,k}$ is an irreducible $\uqsl$-module. 

(ii) If $(\omega,k)\neq (\omega',k')$, the modules $T_{\omega,k}$ and $T_{\omega',k'}$ are not equivalent. 

(iii) Any irreducible finite-dimensional module is equivalent to some $T_{\omega,k}$. 

(iv) Any finite-dimensional $\uqsl$-module is completely reducible. \qed
\end{theorem}

We see that this result is similar to the analogous result for $\frsl(2,\bbC)$, except for the fact that in every possible dimension there are two irreducible modules, distinguished by $\omega$. The modules $T_{1,k}$ are said to be of type I, while modules $T_{-1,k}$ are said to be of type II. The proofs are very similar to the classical proofs.

\begin{theorem}
\label{hd nonroot}
The Dirac cohomology of the irreducible finite-dimensional $\uqsl$-module $(T_{\omega,k},V_k)$ is
\[
H_D(T_{\omega,k})=\bbC v_{-k}\otimes s_{-1} \oplus \bbC v_k\otimes s_1.
\]
The $\uqk$-action on $H_D(T_{\omega,k})$ is given by
\[
K (v_{-k}\otimes s_{-1})=\omega q^{-2k-1}\, v_{-k}\otimes s_{-1};\qquad K (v_k\otimes s_1)=\omega q^{2k+1}\, v_k\otimes s_1.
\]
Note that this action of $K$ is in fact the action of $\delta(K)$, as the algebra $\uqk$ acts on $H_D(M)$ through $\delta(K)$.
\end{theorem}
\pf
This follows from (\ref{Diraccoh}). By (\ref{def vlo}), it is clear that $\ker(E)$ is spanned by $v_k$, and that $v_k\notin \im(F)$. Similarly, $\ker(F)$ is spanned by $v_{-k}\notin\im(E)$. The claim follows.
\epf

\begin{remark}
{\rm
It is instructive to check Theorem \ref{vogconjmod} for the modules $T_{\omega,k}$. One can compute the action of $\Cas_q$ on any $v_m\in V_k$ directly from \eqref{def vlo}, and obtain the scalar
\eq
\label{cas scalar}
\frac{2}{(q-\qi)^2}(\omega(q^{2k+1}+q^{-2k-1})-(q+\qi)).
\eeq
On the other hand, one can write
\[
\zeta(\Cas_q)=\frac{2}{(q-\qi)^2}((\delta(K)+\delta(\ki))-(q+\qi)),
\]
and by Theorem \ref{hd nonroot} this acts on $v_k\otimes s_1$ or on $v_{-k}\otimes s_{-1}$ by the same scalar \eqref{cas scalar}.
}
\end{remark}

%%%%%%%%%%%%%%%%%%%%%%%%%%%%%%%%%%%%%%%%%%%%%%%%%%%%%%%%%%%%%%%%%
\section{Irreducible $\uqsl$-modules for $q$ a root of unity}\label{modrootof1}\label{section5}

Let us now assume $q$ is a primitive $p'$-th root of unity, and set $p=p'$ if $p'$ is odd, $p=p'/2$ if $p$ is even. We are again following \cite{KS}.

By Proposition \eqref{center}, the center ${\mathfrak z}$ 
of $U_q(\mathfrak{sl}_2)$ is generated by $E^p, F^p, K^p, K^{-p}$ and $\Cas_q$. The Poincar\'e-Birkhoff-Witt Theorem then implies that any irreducible $\uqsl$-module is 
finite-dimensional, of dimension at most $p$. 

We recall that the two sided Hopf ideal $J$ in $\uqsl$,
generated by $\langle E^p, F^p, K^p-1 \rangle$, gives rise to the finite-dimensional Hopf algebra called the reduced, or small, quantum group $\reduqsl=\uqsl/J$.

To describe the classification of irreducible $\uqsl$-modules in this case, we retain the notation $(T,V)$ for the representation $T$ on the complex vector space $V$. 
Since $E^p$ and $F^p$ are central elements, $T(E^p)$ and $T(F^p)$ act on $V$ by scalar multiples 
of the identity endomorphism and the classification list splits into four cases, according to these scalars being zero or not.

We start with the representations $(T_{\omega,k},V_k)$ defined in exactly the same way as for $q$ not a root of unity. So $\omega\in\{1,-1\}$, $k\in
\frac{1}{2}\bbN_0$,  
$V_{k}=\langle v_{-k},\ldots ,v_k \rangle$ is a $(2k+1)$-dimensional complex vector space, and the action is given by equations (\ref{def vlo}). (Recall that $v_{k+1}=v_{-k-1}=0$.)

Another, family of representations consists of representations $(T_{a,b,\lambda},V)$, $a,b,\lambda\in{\mathbb C}$, $\lambda\neq 0$. Each of these representations is defined on 
a $p$-dimensional complex vector space $V=\langle v_0,v_1,\ldots ,v_{p-1}\rangle$, and the action is given by the following equations:
\begin{align}\label{Tabl}
& T_{a,b,\lambda}(K)v_m=\lambda q^{-2m}v_m,\qquad\ m\in\{0,\ldots ,p-1\}; \nonumber \\
& T_{a,b,\lambda}(E)v_0=av_{p-1}; \nonumber \\
& T_{a,b,\lambda}(E)v_m=
\big(ab+[m]_q\frac{\lambda q^{1-m}-\lambda^{-1}q^{m-1}}{q-q^{-1}}\big)v_{m-1},
\qquad m>0; \nonumber \\
& T_{a,b,\lambda}(F)v_{p-1}=bv_{0}; \nonumber \\ 
& T_{a,b,\lambda}(F)v_m=v_{m+1}, \qquad m<p-1. 
\end{align} 

The following theorem gives a classification of all irreducible $\uqsl$-modules. (As remarked earlier, all of them are finite-dimensional, since $q$ is a root of unity.)

\begin{theorem}
\label{class root}
\begin{enumerate}
\item The representation $T_{\omega,k}$ is irreducible if and only if $2k<p$. It satisfies $T_{\omega,k}(E^k)=T_{\omega,k}(F^k)=0$. 

The representation $T_{0,0,\lambda}$ is irreducible if and only if $\lambda\not=\pm q^m$ for $m\in\{0,1,\ldots, p-2\}$.
It satisfies $T_{0,0,\lambda}(E^k)=T_{0,0,\lambda}(F^k)=0$. 

The representation $T_{a,b,\lambda}$, $a\neq 0,b\neq 0$, is irreducible if and only if 
\[
ab+[m]\frac{\lambda q^{1-m}-\lambda^{-1}q^{m-1}}{q-q^{-1}}\neq 0
\]
for all $m\in\{0,1,\ldots ,p-1\}$. These representations satisfy 
\[
T_{a,b,\lambda}(E^p)\not=0\not=T_{a,b,\lambda}(F^p);
\] 
such representations are called cyclic and they 
have neither a highest weight vector nor a lowest weight vector. 

The representation $T_{0,b,\lambda}$, $b\neq 0$, is  irreducible if and only if $\lambda^p\not= \pm 1$ (recall that also $\lambda\neq 0$). These representations satisfy
\[
T_{0,b,\lambda}(E^p)=0\not=T_{0,b,\lambda}(F^p).
\]
Such representations are called semicyclic and they 
have a highest weight vector, but no lowest weight vectors.

The representation $T_{a,0,\lambda}$, $a\neq 0$, is  irreducible if and only if $\lambda^p\not= \pm 1$. These representations satisfy
\[
T_{a,0,\lambda}(E^p)\neq 0=T_{a,0,\lambda}(F^p).
\]
Such representations are again called semicyclic; they 
have a lowest weight vector, but no highest weight vectors.

\item If $b'=b$, and if 
\[
a'=a+b^{-1}[m]\frac{\lambda q^{1-m}
-\lambda^{-1}q^{m-1}}{q-q^{-1}} \quad\text{ and }\quad \lambda'=q^{-2m}\lambda
\]
for some $m\in\{0,1,2,\ldots, p-1\}$, then $T_{a,b,\lambda}$ is equivalent to $T_{a',b',\lambda'}$.

If $\lambda=\omega q^{p-1}$, then $T_{0,0,\lambda}$ is equivalent to $T_{\omega,\frac{p-1}{2}}$.

There are no other equivalences between any two of the  representations $T_{\omega,k}$ or $T_{a,b,\lambda}$.

\item Any irreducible $\uqsl$-module is equivalent either to some $T_{\omega,k}$ or to some $T_{a,b,\lambda}$.

\end{enumerate}
\qed
\end{theorem}

We note that the finite-dimensional representations 
$T_{0,0,\lambda}$ for $\lambda=\omega q^n$ ($\omega\in\{1,-1\}$ and $n\in\{0,1,\ldots ,p-2\}$), are indecomposable but not irreducible; we will show a concrete example in Section \ref{section6}. In particular, these modules are not completely reducible. So Theorem \ref{cl fd nonroot}(iv) fails when $q$ is a root of unity.

We also note that the irreducible representations
of $\reduqsl$ come from the first class of $\uqsl$-modules
with $T(E^p)=0=T(F^p)$. These representations were first classified by Lusztig \cite{L1}.

\begin{theorem}
\label{hd root}
\begin{enumerate}
\item The Dirac cohomology of the irreducible $\uqsl$-module $T_{\omega,k}$ is
\[
H_D(T_{\omega,k})=\bbC v_{-k}\otimes s_{-1} \oplus \bbC v_k\otimes s_1.
\]
The $\uqk$-action on $H_D(T_{\omega,k})$ is given by
\[
K(v_{-k}\otimes s_{-1})=\omega q^{-2k-1}\, v_{-k}\otimes s_{-1};\qquad K(v_k\otimes s_1)=\omega q^{2k+1}\, v_k\otimes s_1.
\]
\item
The Dirac cohomology of the irreducible $\uqsl$-module $T_{0,0,\lambda}$ is
\[
H_D(T_{0,0,\lambda})=\bbC v_{p-1}\otimes s_{-1}\, \oplus\, \bbC v_0\otimes s_1.
\]
The $\uqk$-action on $H_D(T_{0,0,\lambda})$ is given by
\[
K(v_{p-1}\otimes s_{-1})=q\, v_{p-1}\otimes s_{-1};\qquad K(v_0\otimes s_1)=q\, v_0\otimes s_1.
\]
\item
The Dirac cohomology of the irreducible $\uqsl$-module $T_{a,b,\lambda}$, $a\neq 0$, $b\neq 0$, is zero.
\item
The Dirac cohomology of the irreducible $\uqsl$-module $T_{0,b,\lambda}$, $b\neq 0$, is 
\[
H_D(T_{0,b,\lambda})=\bbC v_0\otimes s_1.
\]
The $\uqk$-action on $H_D(T_{0,b,\lambda})$ is given by
\[
K(v_0\otimes s_1)=q\, v_0\otimes s_1.
\]
\item
The Dirac cohomology of the irreducible $\uqsl$-module $T_{a,0,\lambda}$, $a\neq 0$, is
\[
H_D(T_{a,0,\lambda})=\bbC v_{p-1}\otimes s_{-1}.
\]
The $\uqk$-action on $H_D(T_{a,0,\lambda})$ is given by
\[
K (v_{p-1}\otimes s_{-1})=q\, v_{p-1}\otimes s_{-1}.
\]
\end{enumerate}
\end{theorem}
\pf
As in the case when $q$ is not a root of unity, this is a direct consequence of \eqref{Diraccoh}; the only thing to do is to identify the kernels of $E$ and $F$ and this is clear from the description of the action given by (\ref{def vlo}) and (\ref{Tabl}).
\epf

%%%%%%%%%%%%%%%%%%%%%%%%%%%%%%%%%%%%%%%%%%%%%%%%%%%%%%%%%%%%%%%
\section{Further examples}\label{section6}

In this section we discuss a few other examples of the 
computation of Dirac cohomology, in cases which are beyond the scope of finite-dimensional irreducible $\uqsl$-modules in Section \ref{modnrootof1} and Section 
\ref{modrootof1}. We present an example of a $\uqsl$-module with infinite-dimensional Dirac cohomology, for $q$ a primitive third root of unity.
We also determine the Dirac cohomology of
a finite-dimensional reducible indecomposable module.

\begin{example}
{\rm
The irreducible finite-dimensional representations we have described in Section \ref{section4}
and Section \ref{section5} are objects of the BGG category ${\fam2 O}_q$ of 
$\uqsl$-modules, cf. \cite{AnMa}. Some examples of infinite-dimensional objects 
in ${\fam2 O}_q$ are the standard (or Verma of type I) $U_q({\mathfrak sl}_2)$-modules 
$V_\lambda$, $\la\in\bbC$. Each $V_\lambda$ is a complex vector space 
$V_\lambda=\bigoplus\limits_{m\in{\mathbb N}_0}\bbC  v_{\lambda-2m}$
equipped with the representation $T_\lambda$ of 
$U_q(\mathfrak{sl}_2)$ given by  
\begin{align}\label{vermarepr}
& T_\lambda(K^{\pm 1})v_{\la-2m}=q^{\pm(\lambda-2m)}v_{\la-2m},\nonumber \\
& T_\lambda(F)v_{\la-2m}=[m+1]v_{\la-2m-2}, \nonumber \\
& T_\lambda(E)v_{\la-2m}=[\lambda-m+1]v_{\la-2m+2}
\end{align}
for $m\in{\mathbb N}_0$. If $q$ is not a root of unity, these modules are highest weight 
modules, with highest weight vector $v_\la$ of 
weight $q^\lambda$, annihilated by $T_\lambda(E)$.

The irreducible finite-dimensional modules of Section \ref{section4} can be obtained as quotients of Verma modules for integral $\la$; we have however changed the indexing of the basis and the constants in the action, so that this is not immediately obvious. (The indexing and the constants we are using now resemble more the usual conventions for classical $\frsl(2,\bbC)$-modules.)

Using \eqref{Diraccoh}, it is easy to see that 
\[
H_D(V_\lambda)=\bbC v_\lambda\otimes s_1,
\]
with the $\uqk$-action given by
\[
K(v_\lambda\otimes s_1) = q^{\la+1} v_\lambda\otimes s_1.
\]

If $q$ is not a root of unity, other infinite-dimensional modules in ${\fam2 O}_q$, e.g., the projective modules $P_\lambda$, can be constructed and their 
Dirac cohomology computed as in the classical case for $U({\mathfrak sl}_2)$. For such modules it makes sense to introduce the notion of higher Dirac cohomology, see \cite{PS}.
}
\end{example}

\begin{example}
{\rm
We now specialize the Verma modules of the previous example to the case when $q$ is a root of unity. The resulting 
representation is no longer generated by a cyclic vector. 
\begin{enumerate}
\item
Let $\lambda=0$ and let $q$ be a third primitive root of unity, $q^3=1$.
Then 
\begin{align}
& \mathrm{Ker}(E)=\bbC v_0 \oplus 
\bigoplus\limits_{i\in\bbN_0}
\bbC v_{-2-6i}; \nonumber \\
& \mathrm{Ker}(F)=\bigoplus\limits_{i\in\bbN_0}
\bbC v_{-4-6i}; \nonumber \\
&\im(E)= \bigoplus\limits_{i\in\bbN_0}
\bbC v_{-2-6i}\oplus \bigoplus\limits_{i\in\bbN_0}
\bbC v_{-4-6i};\nonumber \\
& \im(F)=\bigoplus\limits_{i\in\bbN_0}
\bbC v_{-2-6i}\oplus \bigoplus\limits_{i\in\bbN_0}
\bbC v_{-4-6i} =\im(E).
\end{align}
By \eqref{Diraccoh}, 
\begin{align}
H_D(V_0)= {\mathbb C}v_0 \otimes s_1,
\end{align} 
with $\uqk$-action given by
\[
K(v_0 \otimes s_1) = q v_0 \otimes s_1.
\]

\item
Let $\lambda=1$ and let $q$ be a third primitive root of unity, $q^3=1$. Analogously to the previous case, 
there is a $\uqk$-module isomorphism
\begin{align}
H_D(V_1)\cong {\mathbb C}v_1 \otimes s_1.
\end{align}
\item
Let $\lambda=2$ and let $q$ be a third primitive root of unity, $q^3=1$. Then there is a $\uqk$-module isomorphism
\begin{align}
H_D(V_2)\cong \bigoplus\limits_{i\in{\mathbb N}_0}
{\mathbb C}v_{2-6i}\otimes s_1\,\,\oplus\,\, 
\bigoplus\limits_{i\in{\mathbb N}_0}
{\mathbb C}v_{-2-6i}\otimes s_{-1} .
\end{align}
\end{enumerate}
}
\end{example}
\begin{example}
{\rm
We compute the Dirac cohomology of a finite-dimensional reducible indecomposable $U_q({\mathfrak sl}_2)$-module. Let $q$ be a primitive third root of unity, $q^3=1$. Let $V=\langle v_0,v_1,v_2 \rangle$ 
be a complex vector space of dimension $3$ and $(T_{0,0,1},V)$ be one of the representations
described by \eqref{Tabl}:
\begin{align}
& T_{0,0,1}(K)v_m=q^{-2m}v_m,\qquad\quad m\in\{0,1,2\};
\nonumber \\
& T_{0,0,1}(F)v_0=v_1,\quad T_{0,0,1}(F)v_1=v_2,\quad T_{0,0,1}(F)v_2=0;
\nonumber \\
& T_{0,0,1}(E)v_0=0,\quad T_{0,0,1}(E)v_1=0,\quad T_{0,0,1}(E)v_2=-(q+q^{-1})v_1.
\end{align}
Note that since $q$ is a third primitive root of unity, $q+q^{-1}=-1\not=0$. We see that the subspace $W=\langle v_1,v_2\rangle \subset V$ is a submodule and 
the subspace $\bbC v_0 \subset V$ can be identified with the  quotient module $V/W$. Using \eqref{Diraccoh} again, we conclude
\begin{align}
H_D(T_{0,0,1})= 
\bbC v_2\otimes s_{-1}\,\,\oplus\,\,
\bbC v_0\otimes s_{1}.
\end{align}
}
\end{example}

%%%%%%%%%%%%%%%%%%%%%%%%%%%%%%%%%%%%%%%%%%%%%%%%%%%%%%%%%%%%%%%%%%%%%%%%%

\end{document}